\documentclass{siamart1116}

\usepackage{amsmath}
\usepackage{amsfonts}
\usepackage{amssymb}
\usepackage{url}
\usepackage{latexsym}
\usepackage{graphicx}
\usepackage{enumerate}
\usepackage[shortlabels]{enumitem}
\usepackage{mathtools,mathdots}
\usepackage{hyperref}

\usepackage{amsmath, amssymb,latexsym}
\usepackage{color}
\usepackage{xcolor}
\usepackage{fullpage}

\def\F{\mathbb{F}}
\def\G{\mathbb{G}}

\def\C{\mathbb{C}}

\def\Z{\mathbb{Z}}
\def\Q{\mathbb{Q}}
\def\CS{\mathrm{cs}}
\def\nnull{\mathrm{null}}
\def\A{\mathcal{A}(\Omega)}

\DeclareMathOperator{\rank}{rank}

\DeclareMathOperator{\adj}{adj}
\DeclareMathOperator{\spann}{span}

\newtheorem{remark}[theorem]{{\sc Remark}}

\newtheorem{example}[theorem]{Example}

\title{Invertible bases and root vectors for\\ analytic matrix-valued functions}
\date{}
\author{Vanni Noferini\footnote{Department of Mathematics and Systems Analysis, Aalto University,  P.O. Box 11100, FI-00076, Aalto, Finland, \texttt{vanni.noferini@aalto.fi}. Supported by an Academy of Finland grant (Suomen Akatemian p\"{a}\"{a}t\"{o}s 331240).}  }

\begin{document}

\maketitle

% REQUIRED
\begin{abstract}
We revisit the concept of a minimal basis through the lens of the theory of modules over a commutative ring $R$. We first review the conditions for the existence of a basis for submodules of $R^n$ where $R$ is a B\'{e}zout domain. Then, we define the concept of invertible basis of a submodule of $R^n$ and, when $R$ is an elementary divisor domain, we link it to the Main Theorem of [G. D. Forney Jr., SIAM J. Control 13, 493--520, 1975]. Over an elementary divisor domain, the submodules admitting an invertible basis are precisely the free pure submodules of $R^n$. As an application, we let $\Omega \subseteq \C$ be either a connected compact set or a connected open set, and we specialize to $R=\A$, the ring of functions that are analytic on $\Omega$. We show that, for any matrix $A(z) \in \A^{m \times n}$, $\ker A(z) \cap \A^n$ is a free $\A$-module and admits an invertible basis, or equivalently a basis that is full rank upon evaluation at any $\lambda \in \Omega$. Finally, given $\lambda \in \Omega$, we use invertible bases to define and study maximal sets of root vectors at $\lambda$ for $A(z)$. This in particular allows us to define eigenvectors also for analytic matrices that do not have full column rank.
\end{abstract}

% REQUIRED
\begin{keywords}
Analytic function, eigenvector, free module, minimal basis, invertible basis, root vector, maximal set, pure submodule
\end{keywords}

% REQUIRED
\begin{AMS}
15B33, 15A18, 65H17, 15B99
\end{AMS}

\maketitle

\section{Introduction}

In the context of the theory of polynomial matrices \cite{GLR}, the concept of a \emph{minimal basis} was introduced by D. Forney Jr. in \cite{For75}. Let $\F$ be a field, and denote by $\F[x]$ (resp. $\F(x)$) the ring of univariate polynomials over $\F$ (resp. the field of fractions of $\F[x]$). If $V$ is any subspace of $\F(x)^n$, it is clear that a polynomial basis for $V$ exists. Forney gave the following definition.
\begin{definition}\label{def:minimal}
Let $\mathcal{B}=\{b_1(x),\dots,b_p(x)\}$ be a polynomial basis for $V \subseteq \F(x)^n$. Then, the quantity $\Omega=\sum_{i=1}^p \deg b_i(x)$ is called the order of $\mathcal{B}$. If $\mathcal{B}$ has minimal order among all possible polynomial bases of $V$, it is called a \emph{minimal basis} of $V$.
\end{definition}
Then, Forney showed that the degrees of a minimal basis, called the \emph{minimal indices} of $V$, are uniquely determined by $V$. Moreover, he gave some equivalent characterization of a minimal basis in \cite[Main Theorem]{For75}. To state it, we must define the high coefficient matrix $\{P\}$ of a polynomial matrix $P(x)$ as the matrix whose columns retain the coefficients of the highest order monomials in each column of $P(x)$. For example, over $\F=\Q$,

\[ \left\{ \begin{bmatrix}
1 & 2x^2\\
0 & x+1\\
-1 & x^2+x+1
\end{bmatrix} \right\} = \begin{bmatrix}
1&2\\
0&0\\
-1&1
\end{bmatrix}.\]

\begin{theorem}[Forney]\label{thm:forney}
Let $M(x) \in \F[x]^{n \times p}$. The following are equivalent:
\begin{enumerate}
\item The set of the columns of $M(x)$ is a minimal basis for the subspace they span;
\item (a) For every prime polynomial $p(x) \in \F[x]$, $M(x)$ is nonsingular modulo $p(x)$ and (b) $\{M\}$ is nonsingular;
\item (a) The GCD of all $p \times p$ minors of $M(x)$ is $1$ {and} (b) The largest degree of all the $p \times p$ minors of $M(x)$ is the sum of the degrees of the columns of {$M(x)$};
\item If $b(x)=M(x)a(x)$ is a polynomial vector then (a) $a(x)$ must be polynomial {and} (b) The degree of $b(x)$ is the maximum, for all $1 \leq i \leq p$, of the sum of the degree of the $i$-th column of $M(x)$ and the degree of the $i$-th coefficient of $a(x)$.
\end{enumerate}
\end{theorem}
(For conciseness, we have omitted a fifth equivalent condition from Forney's original theorem, that relates, for all $d$, the minimal indices of $V$ with the dimension over $\F$ of $V_d$, the $\F$-vector space of all polynomial $n$-uples in $V$ of degree $\leq d$.)

Two more remarks on minimal bases are in order. Our first remark is that Forney observed already in \cite{For75} that properties (a) in his Main Theorem are equivalent to the fact that the $\F[x]$-module $M=V \cap \F[x]^n$ is free and that a minimal basis is a basis for $M$. However, Forney mentioned this (correct) fact in passing, without explaining it in detail; a more thorough treatement can be found, e.g., in \cite[Section 2.4]{fh}. Moreover, Forney showed that in Theorem \ref{thm:forney} properties (a) are equivalent to each other (and to the fact that a minimal basis of $V$ is a basis of $M$) and properties (b) are equivalent to each other.

The second remark is that root polynomials and root vectors at $\lambda \in \C$ for a (possibly not full column rank) polynomial matrix $P(x)$ were defined in \cite{DN,Nof12}, extending the original notion for regular polynomial matrices \cite{GLR} by building on the concept of the vector space $\ker_\lambda P(x)$; the latter was in turn defined as the span of a minimal basis of $\ker P(x)$, evaluated at $x=\lambda$. A so-called maximal set of root polynomials is a useful tool as it provides all the local positive partial multiplicities at $\lambda$: see \cite{DN,Nof12} as well as \cite{NV} for an extension to rational matrices. The idea of root vectors has proved useful for theoretical analysis \cite{DN,DN2,Nof12,NV} as well as to analyze conditioning \cite{LN,NNPQ} and to design practical algorithms for singular eigenvalue problems \cite{HMP,KG}. However, a careful inspections of the arguments in \cite{DN,Nof12,NV} shows that to define $\ker_\lambda P(x)$ it suffices in fact to start from a, not necessarily minimal, basis satisfying properties (a) in Forney's Main Theorem \ref{thm:forney}. In other words, properties (b) in Theorem \ref{thm:forney} are not necessary to correctly define root polynomials at every finite point. In this paper, we propose to call a polynomial basis that satifies properties (a) in Forney's Main Theorem \ref{thm:forney} \emph{invertible}, since a matrix whose columns are the vectors of such basis is left invertible (over $\F[x]$). Equivalent characterizations of left invertibility are that the matrix has trivial Smith form, i.e., all the invariant factors are units, or that it is a submatrix of a unimodular matrix; see also \cite[Theorem 3.3]{Zaballa}. 

The main goals of this paper are (1) To {observe that, when $R$ is an elementary divisor domain (but not necessarily a principal ideal domain) the notion of invertible basis is generalized by the bases of certain submodules of $R^n$ known as pure submodules} and (2) To use invertible bases to define a maximal set of root vectors for matrices over the ring of analytic functions on $\Omega$ where $\Omega$ is a connected (open or compact) set. We expect that the resulting theory may be useful, for example, for applications to the nonlinear eigenvalue problem \cite{GT}. We will take a module-theoretical approach, while keeping the exposition accessible for an audience with expertise in numerical linear algebra. In particular, in Section \ref{sec:modules} we explicitly recall all the needed algebraic concepts about modules. In Section \ref{sec:BD}, we { recall several known results concerning submodules of $R^n$ that have useful consequences.}
%show that when $R$ is a B\'{e}zout domain with field of fractions $\G$, then a submodule of $R^n$ is free if and only if it is finitely generated (Theorem \ref{thm:submodulesofBn}).
 If in particular $R$ is an elementary divisor domain, and for any suspace $V \subseteq \G^n$, then $V \cap R^n$ is free and it has an invertible basis (Proposition \ref{cor:minimalatfiniteBn}). This  result is certainly not new: see for example \cite[Theorem 1.13.3]{Friedland}; however, our presentation and in particular the link to properties (a) in Theorem \ref{thm:forney} (see Theorem \ref{thm:plc}) seem to not be well known in the linear algebra community. {If a submodule of $R^n$  admits an invertible basis, then all its bases are also invertible (see Remark \ref{rem:referee1}). Over elementary divisor domains, this property characterizes the pure free submodules of $R^n$ \cite{carmsa,faith,kaplansky1}.} Finally, in Section \ref{sec:maximalsets} we specialize to the ring $\A$ of analytic functions over a connected, and either open or compact, set $\Omega$, and we use the existence of invertible basis to define and characterize maximal sets of root vectors for matrices over $\A$. This allows us in particular to define eigenvectors even for analytic matrices that do not have full column rank: in that case, an eigenvector can be seen as a certain equivalence class.

\section{Commutative rings, modules, and free modules}\label{sec:modules}
In this section we refer to some basic algebra material: the notions that we recall, and more, can be found for example in \cite{clark,faith,Friedland,fh,hs,Kaplansky}. Throughout, we consider a commutative ring $R$ with unity, and we assume $1 \neq 0$ so that $R$ is not the trivial ring $\{0\}$, because it contains at least two distinct elements. A \emph{module} over the commutative ring $R$, or $R$-module, is a non-empty set $M$ which is a commutative group with respect to the addition $+$ (in particular, there is a zero element $0 \in M$ and every $x \in M$ has an inverse $-x \in M$ such that $x+(-x)=0$) and is endowed with a scalar multiplication $\cdot: R \times M \rightarrow M$ that satisfies the following axioms: (A1) 
$a \cdot (x+y) = a \cdot x + a \cdot y \ \mathrm{for} \ \mathrm{all} \ x,y \in M, a \in R;$ (A2)
$(a+b) \cdot x = a \cdot x + b \cdot x \ \mathrm{for} \ \mathrm{all} \ x \in M, a,b \in R;$ (A3)
$(ab) \cdot x = a \cdot (b \cdot x) \ \mathrm{for} \ \mathrm{all} \ x \in M, a,b \in R;$ (A4) $1 \cdot x=x \ \mathrm{for} \ \mathrm{all} \ x \in M.$

More informally, axioms (A1)--(A4) say that a module is to a commutative ring what a vector space is to a field. (In fact, a vector space is an $R$-module in the case of $R$ being a field.) In particular, the concepts of linear independence and of linear span both still make sense. Namely, we say that $x_1,\dots,x_n \in M$ are \emph{linearly independent} if, given $a_1,\dots,a_n \in R$, $\sum_{i=1}^n a_i \cdot x_i = 0 \Rightarrow a_i=0 \ \forall \ i$. Moreover, we say that the module $M$ is (finitely) generated by the elements $y_1,\dots,y_m \in M$ (called \emph{generators}) if, for every $x \in M$, there exist $b_1,\dots,b_m \in R$ such that $x=\sum_{i=1}^m b_i \cdot y_i$. Linear independence can be applied also to infinite sets of elements of $M$: in this case, one requires that every possible finite subset is linearly independent. Similarly, a module $M$ may also be infinitely generated: in this case, it is necessary to have an infinite number of generators but every element of $M$ is still required to be an $R$-linear combination of finitely many generators. If $S \subset M$ is a set of generators of the $R$-module $M$, we write $M=\mathrm{span}_R S$.

However, when combining these two concepts, a striking difference between modules and vectors spaces emerges. Indeed, define a \emph{basis} of a module as a linearly independent set of generators. Every vector space has a basis\footnote{Of course, in the infinite dimensional case, this result is equivalent to the axiom of choice.} but the same is not true for modules.
\begin{example}
Take $R=\mathbb{Z}$. Then, it is easy to verify that $M=\mathbb{R}$ is a $\mathbb{Z}$-module. Assume for a contradiction that it has a basis $\mathcal{B}$ and let $b_0 \in \mathcal{B} \subseteq \mathbb{R}$. Suppose that $2 \leq k \in \mathbb{Z}$, then by assumption there exist $b_1,\dots,b_n \in \mathcal{B}$ such that, for some $z_0,z_1,\dots,z_n \in \mathbb{Z}$, 
\[ \frac{b_0}{k} = \sum_{i=0}^n z_i b_i \Rightarrow b_0 - \sum_{i=0}^n (k z_i) b_i = 0,\]
contradicting the linear independence of the basis. (To deal with the case where $z_0 \neq 0$, note that $1-kz_0 \neq 0$ as otherwise $|k|=1$, contradicting $k \geq 2$.)
\end{example}
Nevertheless, some modules do have a basis. A module that has a basis is called \emph{free}. Moreover, for convenience of exposition, we will conventionally agree that the trivial module $M=\{0\}$ is also free, with its basis being the empty set; this excentricity actually makes sense if we formally agree that the sum of an empty set of elements of $M$ yields $0$, and in the following it will allow us to avoid the clumsy need to repetitively exclude the trivial module when making formal statements about free modules.
\begin{proposition}\label{prop:Rnisfree}
Let $R$ be a commutative ring. $R^n$ (the set of all possible $n$-uples of elements of $R$) is a free $R$-module.
\end{proposition}
\begin{proof}
It is clear that the canonical basis $\{e_i\}_{i=1}^n$, where $(e_i)_j=1$ if $i=j$ and $(e_i)_j=0$ if $i \neq j$, is a basis of $R^n$.
\end{proof} 

The example of Proposition \ref{prop:Rnisfree} is in fact unique up to isomorphism, in the sense that every finitely generated free $R$-module with a basis of cardinality $n$ is isomorphic to $R^n$.

We conclude this section by recalling some further nomenclature and basic properties \cite{clark,faith,Friedland,fh,hs,Kaplansky}. In the special case $M \subseteq R$, the $R$-module $M$ is called an \emph{ideal}.  Recall also that a commutative ring $R$ is called an \emph{integral domain} if, for all $a,b \in R$, $ab=0$ and $b\neq 0$ imply $a=0$. An integral domain is called a \emph{principal ideal domain} (PID) if every ideal of $R$ is principal, i.e., generated by a single element of $R$;  and it is called a \emph{B\'{e}zout domain} (BD) if every finitely generated ideal of $R$ is principal. Obviously, every PID is a BD by definition, but the converse is not true. For example, the ring $\mathcal{A}(\C)$ of entire functions is a BD \cite{Friedland} but not a PID. To see that $\mathcal{A}(\C)$ is not a PID, define $J \subset \mathcal{A}(\C)$ as the set of entire functions $f(z)$ that are zero at all but finitely many Gaussian integers: it is easy to verify that $J$ is an ideal which is not finitely generated (and in particular not principal). Finally, an integral domain $R$ is called an \emph{elementary divisor domain} (EDD) if the following theorem holds.
\begin{theorem}[Smith]\label{thm:smith}
Let $R$ be an elementary divisor domain and $A \in R^{m \times n}$. Then there exist two unimodular (that is, invertible over $R$) matrices $U \in R^{m \times m}$ and $V \in R^{n \times n}$ such that $A=USV$ where $S \in R^{m \times n}$ is diagonal and such that the $(i,i)$ element of $S$ divides the $(i+1,i+1)$ element of $S$, for all $i < \min\{m,n\}$.
\end{theorem}
{  A matrix $S$ having the properties stated in Theorem \ref{thm:smith} is called a \emph{Smith form} of $A$, and it is uniquely determined by $A$ up to multiplication of each diagonal element by a unit of $R$. The diagonal elements of $A$ are called \emph{invariant factors} of $A$, and the product of the first $k$ invariant factors is called the $k$th \emph{determinantal divisor} of $A$ and it is equal to the GCD of all the $k \times k$ minors of $A$. Both invariant factors and determinantal divisors are uniquely determined by $A$ up to multiplication by a unit of $R$.}

It can be proved that every PID is an EDD and that every EDD is a BD \cite{Friedland}, while to our knowledge it is still an open problem to decide whether there is a BD which is not an EDD or the definitions of EDD and BD are equivalent \cite{Lorenzini}.
\section{Submodules of $R^n$ when $R$ is a PID, an elementary divisor domain, or a B\'{e}zout domain}\label{sec:BD}
In this section, we study submodules of $R^n$ where $R$ is a commutative ring; we will often make further assumptions on $R$ such as being a PID, EDD or BD. Most of the results in this section can be traced elsewhere, e.g., in \cite{brown,faith,For75,Friedland,fh,hs}, although our derivation does not strictly follow those sources.

Let $R$ be a commutative ring and $A \in R^{m \times n}$. The range, or column space, of $A$ is defined as $\CS(A):=\{ x \in R^m : \exists \ y \in R^n \ s.t. \ x=Ay\}$. It is a simple exercise to verify that $\CS(A) \subseteq R^m$ is an $R$-module. The null space (over $R$) of $A$ is denoted by $\mathrm{null}(A):=\{y \in R^n : Ay=0\}$; when $R$ is an integral domain we notationally distinguish $\mathrm{null}(A)$ from $\ker A$, that in this paper denotes the null space, or kernel, of $A$ over $\G$, the field of fractions of $R$. (For example, if $R=\F[x]$ is the ring of polyomials then $\G=\F(x)$ is the field of rational functions.) It is clear that, when $R$ is an integral domain, $\mathrm{null}(A)=\ker A \cap R^n$.

 The following result \cite[Theorem 5.10]{brown} is very general as no assumption on $R$, beyond being a commutative ring, is needed.
\begin{theorem}[Brown]\label{thm:brown}
Let $R$ be a commutative ring and let $M$ be a finitely generated $R$-module. Suppose that $\{m_1,\dots,m_k\} \subseteq M$ is a linearly independent set and $\{p_1,\dots,p_n\} \subseteq M$ is a set of generators of $M$. Then $k \leq n$. Moreover, if $k=n$, then $M$ is free and $\{p_1,\dots,p_n\}$ is a basis of $M$.
\end{theorem}
An immediate corollary of Theorem \ref{thm:brown} is that, if a matrix $A \in R^{n \times p}$ over a commutative ring has linearly independent (over $R$) columns, then $\CS(A)$ is a free module and the set of the columns of $A$ is a basis. (This also follows from the first isomorphism theorem for modules since in this case $\CS(A) \cong R^p / \nnull(A)  = R^p$.)  Generally, however, neither the range nor the null space of a matrix need to be free modules.
\begin{example}
Let $A=2 \in \Z_4$. Clearly, $\nnull(A)=\{0,2\}$ is generated by $2$ and $2$ is the unique possible generator. However, $\{2\}$ is not a basis for {$\nnull(A)$} since it is not a $\Z_4$-independent set (indeed  in $\Z_4$ it holds $2 \cdot 2 = 0$). Moreover, $\CS(A)=\nnull(A)$ is also not free.
\end{example}
Luckily, though, over principal ideal domains the situation is much simpler.
\begin{theorem}\label{thm:submodule}
Let $R$ be a principal ideal domain. Then, every submodule of $R^n$ is free.
\end{theorem}
\begin{proof}
By Proposition \ref{prop:Rnisfree}, $R^n$ is a free $R$-module. On the other hand, every submodule of a free $R$-module is free if $R$ is a PID \cite[Theorem 5.1]{hs}.
\end{proof}
\begin{corollary}\label{cor:minimalatfinite}
Let $R$ be a principal ideal domain and $\G$ its field of fractions. For every subspace $V \subseteq \G^n$, the set $M=V \cap R^n$ is a free $R$-module.
\end{corollary}
\begin{proof}
In view of Theorem \ref{thm:submodule}, and since $M \subseteq R^n$,  it is enough to show that $M$ is a module. It is clear that $M$ is a group with respect for addition, for $x,y \in M$ implies both $x+y \in V$ and $x+y \in R^n$; and it is also clear that if $x \in M$ and $p \in R$ then $px \in R^n$ and $px \in V$. Finally, properties (A1)--(A4) are also immediate because $M \subseteq R^n$.
\end{proof}
\begin{corollary}\label{cor:mackey}
Let $R$ be a principal ideal domain and let $S$ be a subset of $R^n$. Then, $M=\mathrm{span}_R S$ is a free $R$-module.
\end{corollary}
\begin{proof}
Similarly to Corollary \ref{cor:minimalatfinite}, it is easy to verify that $M$ is a module. Hence, $M$ is a free module by Theorem \ref{thm:submodule}.
\end{proof}

Corollary \ref{cor:minimalatfinite} and Corollary \ref{cor:mackey} provide, for example, a somewhat more natural enviroment for the results related to the filtration approach to minimal bases \cite{M21}: in that context, one can consider more directly a nested sequence of submodules of $\F[x]^n$ rather than of subspaces of $\F(x)^n$. 

We now turn to B\'{e}zout domains, where there are some complications with respect to principal ideal domains. Nevertheless, the complications are not too hard to overcome.

\begin{theorem}\label{thm:submodulesofBn}
Let $R$ be a B\'{e}zout domain and let $M  \subseteq R^n$ be a submodule. Then, $M$ is a free $R$-module if and only if it is finitely generated.
\end{theorem}

\begin{proof}
We first prove necessity. Suppose that $M$ is free and let $\mathcal{B}$ be a basis of $M$. If $\# \mathcal{B} > n$, let the columns of $A \in R^{n \times p}$ be any finite subset of $\mathcal{B}$ of cardinality $p > n$: then it suffices to take any nonzero vector $c \in \nnull(A)$ to immediately show that $\mathcal{B}$ must be linearly dependent. Hence, $\# \mathcal{B} \leq n$, and in particular $\mathcal{B}$ is a finite set thus showing that $M$ is finitely generated.

{The sufficiency is a consequence of \cite[Theorem 1]{kaplansky1}, see also \cite[proof of Proposition 1.8]{nove}.}

\end{proof}

{Among free modules, some have bases that retain, in a generalized sense, properties (a) in Forney's Theorem \ref{thm:forney}. To make this statement more precise, we} now formalize the definition of an invertible basis anticipated in the introduction.

\begin{definition}\label{def:invbasis}
Let $R$ be a commutative ring and $M \subseteq R^n$ be a free $R$-module with basis $\mathcal{B}$. We say that $\mathcal{B}$ is an \emph{invertible basis} if there is a matrix $A \in R^{n \times p}$ such that (1) the columns of $A$ are the elements of $\mathcal{B}$ (2) $A$ is left invertible over $R$, i.e., there exists $L \in R^{p \times n}$ such that $LA=I_p$.
\end{definition}

\begin{remark}\label{rem:referee1}
{ Note that admitting an invertible basis is a property of the module $M$ itself: either every basis of $M$ is invertible or no one is. Indeed, suppose that $A_1,A_2$ are  two matrices whose columns are the elements of two bases of $M$. Then, $A_1 = A_2 Z$ for some unimodular matrix $Z$, and hence $A_1$ is left invertible if and only if $A_2$ is. 

If $R$ is an elementary divisor domain, the modules whose bases are invertible are precisely those that in ring theory are called \emph{pure submodules}: see, e.g.,\cite[p.79]{carmsa}, \cite[p. 70]{faith} or \cite{kaplansky1} for a definition and \cite[p.80]{carmsa} for some properties. Hence, when $R$ is an EDD, Definition \ref{def:invbasis} is equivalent to $M \subseteq R^n$ is a free pure submodule of $R^n$ and $\mathcal{B}$ being a basis of $M$. Moreover, $M$ is a free pure submodule of $R^n$ if and only if every basis of $M$ is invertible.}
\end{remark}

\begin{proposition}\label{cor:minimalatfiniteBn}
Let $R$ be an elementary divisor domain and $\G$ its field of fractions. For every subspace $V \subseteq \G^n$, the set $M=V \cap R^n$ is a free $R$-module. Moreover, {$M$ is a pure submodule of $R^n$, i.e., every basis of $M$ is invertible.}
\end{proposition}
\begin{proof}
That $M$ is a module can be easily proved as in Corollary \ref{cor:minimalatfinite}. Since every EDD is a BD, and in view of Theorem \ref{thm:submodulesofBn}, we must show that $M$ is finitely generated: we will do so constructively, and the procedure will yield an invertible basis as a bonus. {(Note that by Remark \ref{rem:referee1} it suffices to exhibit one.)} Let us start with a basis (over $\G$) of $V$ with elements in $R^n$: this can be easily constructed starting by any basis and scaling each element by {its} least common denominator\footnote{This is possible because every BD is a GCD domain, and it is known \cite[Theorem 2]{K03} that in a GCD domain every pair of elements have both a GCD and a LCM.}. Let $A \in R^{n \times p}$ be the matrix whose columns are the elements of such a basis. Then, $A=QSZ$ where $Q \in R^{n \times n}$ and $Z \in R^{p \times p}$ are unimodular while $S \in R^{n \times p}$ is a Smith form of $A$ (see Theorem \ref{thm:smith}). Letting $Q_p \in R^{n \times p}$ be the matrix containing the leftmost $p$ columns of $Q$ and $S_p \in R^{p \times p}$ be the matrix containing the top $p$ rows of $S$, we also have $A=Q_p S_p Z$. We claim that the set of the columns of $Q_p$ is an invertible basis of $M$. Linear independence is clear, as $Q_p$ is a submatrix of a unimodular matrix. It remains to show that $M$ is generated over $R$ by the columns of $Q_p$. To this goal, note that the columns of $A$ span $V$ over $\G$, and hence for all $v \in M \subseteq V$ there is $c \in \G^p$ such that $v = A c$. Thus, $v=Q_p (S_p Z c)$. But $Q_p$ is left invertible (over $R$) by construction, so by letting $L \in R^{p \times n}$ be a left inverse it holds $S_p Z c = L v \in R^p$: and since $v$ is a generic element of $M$, the columns of $Q_p$ are a set of generators of $M$.
\end{proof}

{Properties (a) in Forney's Theorem \ref{thm:forney} include} property 4(a) that, in the case $R=\F[x]$, was labelled ``polynomial linear combination property" in \cite{M21}. To state Theorem \ref{thm:plc}, recall that in a commutative ring $R$, given $a,b,r \in R$, we write $a \equiv b \mod r$ to mean that $r$ divides $a-b$, i.e., that there is $d \in R$ such that $a-b=rd$. This notation extends elementwise to matrices over $R$: for $A,B \in R^{m \times n}$ we write $A \equiv B \mod r$ if there exists $D \in R^{m \times n}$ such that $A-B=rD$.

\begin{theorem}\label{thm:plc}
Let $R$ be an elementary divisor domain and $\G$ its field of fractions. Let $V \subseteq \G^n$ be a suspace and suppose that set of the columns of the matrix $Q \in R^{n \times p}$ is a basis for $V$. Then the following are equivalent:
\begin{enumerate}
\item The set of the columns of $Q$ is an invertible basis for the module $V \cap R^n$;
\item For all non-units $r \in R$ and for all $Q' \in R^{n \times p}$ such that $Q \equiv Q' \mod r$, $Q'$ has full column rank;
\item The GCD of all $p \times p$ minors in $Q$ is $1$;
\item If $b=Qa \in R^n$ for some $a \in \G^p$, then $a \in R^p$;
\item $Q$ has trivial Smith form, i.e., a Smith form of $Q$ is $S=\begin{bmatrix}
I_p\\
0
\end{bmatrix}$.
\end{enumerate}
\end{theorem}
\begin{proof}
\begin{itemize}
\item[$1 \Rightarrow 2$] By assumption there is $L \in R^{p \times n}$ such that $LQ=I_p$. Let $r \in R$ and let $D,Q' \in R^{p \times n}$ satisfy $rD=Q'-Q$, then $LQ'=I_p + r LD \Rightarrow \det(LQ') \equiv 1 \mod r$. Therefore, if $\det(LQ')=0$ then $1=rd$ for some $d \in R$ implying that $r$ is a unit. Hence, if $r$ is not a unit then $\det(LQ')\neq 0$ and thus $Q'$ has full rank.
\item[$2 \Rightarrow 3$] Write $Q=USV$ with $U,V$ unimodular and $S$ in Smith form. Let $g \in R$ be the GCD of all $p \times p$ minors in $Q$, and suppose $g$ is not a unit. Then, by Theorem \ref{thm:smith}, the $p$-th invariant factor of $Q$, or equivalently the $(p,p)$ element of $S$, cannot be a unit, say, $S_{pp}=r$ for some non-unit $r \in R$.  If $S'=S-re_pe_p^T$ is constructed from $S$ by replacing $S_{pp}$ by $0$, define $Q'=US'V$. Clearly $Q'$ is rank deficient and $Q-Q' = r Ue_pe_p^TV \equiv 0 \mod r$.
\item[$3 \Rightarrow 5$] If the GCD of all $p \times p$ minors of $Q$, that is the $p$-th determinantal divisor, is $1$ then by Theorem \ref{thm:smith} all the invariant factors of $Q$ must be units of $R$, and hence $\begin{bmatrix}
I_p\\
0
\end{bmatrix}$ is a Smith form of $Q$.
\item[$5 \Rightarrow 1$] Since $S^T S=I_p$ and $Q=USV$ for some unimodular $U,V$, a left inverse of $Q$, say $L$, can be constructed as $L=V^{-1}S^TU^{-1}$.
\item[$1 \Rightarrow 4$]
Suppose $Q$ is an invertible basis; then, by definition, there exists $L \in R^{p \times n}$ such that $LQ=I_p$. Hence, $a=Lb \in R^p$. 

\item[$4 \Rightarrow 2$] Let $r \in R$ be a non-unit and suppose that there exists $Q' \in R^{n \times p}$ which is rank deficient and satisfies $Q' \equiv Q \mod r$. Let $c \in \nnull(Q')$ be such that the GCD of the elements of $c$ is $1$. (It is clear that such $c$ exists, by taking any nonzero element of $\nnull(Q')$ and dividing it by the GCD of its entries.) Let $D \in R^{p \times n}$ satisfy $rD=Q-Q'$, then $Qc=Q'c+(rD)c=D(rc)$. Hence, if $a=c/r \in \G^p$, we have that $a \not \in R^p$ but $Qa=Dc \in R^n$.
\end{itemize}
\end{proof}

\begin{remark}\label{rem:referee2}
{ The equivalences $1 \Leftrightarrow 3 \Leftrightarrow 5$ in Theorem \ref{thm:plc} also admit a more direct (and more module-theoretic) proof as they are immediate consequences of the fact that  $V \cap R^n$ is a pure submodule of $R^n$.}
\end{remark}

\begin{remark}An important subtlety is that property 2 in Theorem \ref{thm:plc} is being full rank modulo every non-unit, while property 2(a) in Theorem \ref{thm:forney} is being full rank modulo every prime. Since every prime is a non-unit, it is clear that the former implies the latter. 

Over a PID, such as the ring of univariate polynomials over a field, then every element of the ring admits a unique (up to permutations and multiplications by units) prime decomposition; hence, over a PID the converse implication is also true and the two properties are thus equivalent, because every element of a PID is either a unit or divisible by a prime. For example, Theorem \ref{thm:plc} implies that an integer basis for a subspace $V\subseteq\Q^n$ is invertible if and only if it is full rank modulo $p$, for every prime number $p$. 

In general, however, property 2 in Theorem \ref{thm:plc} is not equivalent to being full rank modulo every prime if $R$ is not a PID. For example let $\mathbb{A} \subset \C$ be the ring of algebraic integers (roots of \emph{monic} polynomials in $\mathbb{Z}[x]$) which is a BD \cite[Theorem 102]{Kaplansky} and an EDD \cite[Theorem 5]{diciassette} so Theorem \ref{thm:plc} holds. On the other hand, $\mathbb{A}$ has no irreducible elements (and hence no prime elements), because the square root of an algebraic integer is an algebraic integer. It follows that $\mathbb{A}$ contains non-units (say, $2 \in \mathbb{A}$ but $\frac12 \not \in \mathbb{A}$ so $2$ is not a unit of $\mathbb{A}$) that are not divisible by any prime elements. We conclude that, for an EDD that is not a PID, property 2 as stated in Theorem \ref{thm:plc}, and hence being an invertible basis, is possibly stronger than (the analogue of) property 2(a) as stated in Theorem \ref{thm:forney}; indeed there are EDDs in which being full rank modulo every prime is \emph{not} equivalent to being an invertible basis. For example $\{2\}$ is not an invertible basis of $\mathbb{A}$ although $2$ is (vacuously) full rank modulo any prime.
\end{remark}

\section{Invertible bases and root vectors over $\A$}\label{sec:maximalsets}

Throughout this section, we fix a connected set $\Omega \subseteq \C$, assuming that $\Omega$ is either compact or open, and we denote by $\A$ the ring of functions that are analytic on $\Omega$. The spectral theory of analytic matrices is relevant in the applications, and in particular in the context of nonlinear eigenvalue problems or in the context of signal processing: see \cite{GT} for a recent survey on nonlinear eigenvalue problems or \cite{Weiss} for some applications in signal processing.

 We recall that, under the stated assumptions, $\A$ is an elementary divisor domain (the Smith theorem holds) and hence a BD: see for example \cite{Friedland}. Thus, the analysis of Section \ref{sec:BD} is relevant. In fact, in the case where $\Omega$ is a connected compact set, even more strongly $\A$ is actually a PID and a Euclidean domain \cite[Lemma 1.3.7]{Friedland}. We consider $\A^{m \times n}$, the set of $m \times n$ matrices over $\A$. By obvious extension of the properties of scalar analytic functions, every element of $\A^{m \times n}$ admits a convergent power series at any point $\lambda \in \Omega$, say, $A(z)=\sum_{j=0}^{\infty} A_j (z-\lambda)^j$, with $(A_j)_j \subset \C^{m \times n}$. A nonzero analytic matrix $A(z)\neq 0$ has a root of order $k$ at $\lambda$ if $A_0=\dots=A_{k-1}=0 \neq A_k$; or equivalently if $A(z) = (z-\lambda)^k B(z)$ for some $B(z) \in \A^{m \times n}$ such that $B(\lambda)\neq 0$.

Let us now fix a matrix $A(z) \in \A^{m \times n}$. The most basic definition is that of an eigenvalue and of the multiplicities associated with it.
\begin{definition}[Eigenvalues and multiplicities]
We say that $\lambda \in \Omega$ is an \emph{eigenvalue} of the matrix $A(z) \in \A^{m \times n}$ if $\mathrm{rank} A(\lambda) < \mathrm{rank} A(z)$. Moreover, the \emph{partial multiplicities} of an eigenvalue of $A(z)$ are the orders of $\lambda$ as a root of the nonzero invariant factors of $A(z)$. The number of the nonzero partial multiplicities of $\lambda$ is called the \emph{geometric multiplicity} of $\lambda$, and the sum of the partial multiplicities is called the \emph{algebraic multiciplicity} of $\lambda$.
\end{definition}

 The {free pure submodule} $M=\nnull(A(z))$ admits an invertible basis by {Proposition \ref{cor:minimalatfiniteBn}}. Let us take one such basis, arrange its elements as the columns of a matrix and denote the latter by $Q(z)$. We first show that the set of the columns of $Q(z)$ is an invertible basis if and only if no $\lambda \in \Omega$ is an eigenvalue of $Q(z)$.

\begin{proposition}
The set of the columns of the matrix $Q(z) \in \A^{n \times p}$ is an invertible basis for the module they span if and only if $Q(\lambda)$ has full column rank for all $\lambda \in \Omega$.
\end{proposition}
\begin{proof}The statement follows from Theorem \ref{thm:plc}, but for completeness let us give an argument more specific to $\A$. Suppose that the set of the columns of $Q(z)$ is an invertible basis for the module $\CS(Q(z))$, then there exists $L(z) \in \A^{p \times n}$ such that $L(z)Q(z)=I_p$, implying $L(\lambda)Q(\lambda)=I_p$ and hence $Q(\lambda)$ has full column rank. Conversely, if $Q(\lambda)$ has full colum rank for all $\lambda \in \Omega$ then, noting that the units of $\A$ are precisely the analytic functions with no zeros in $\Omega$, the Smith form of $Q(z)$ is trivial, and hence $Q(z)$ is left invertible.
\end{proof}

For brevity, and coherently with common practice in the literature on minimal bases, from now on we may write, e.g., ``$Q(z)$ is an invertible basis of the module $M$" as shorthand to mean ``the set of the columns of $Q(z)$ is an invertible basis for the module $M$".

\begin{definition}
Let $A(z) \in \A^{m \times n}$ and suppose that $Q(z) \in \A^{n \times p}$ is an invertible basis of $\nnull(A(z))$. Given $\lambda \in \Omega$, we define $\ker_\lambda A(z):=\mathrm{span} \ Q(\lambda) \subseteq \C^n$.
\end{definition}

It is clear that $\ker_\lambda A(z)$ is well defined, as if $Q(z)$ and $T(z)$ are two invertible bases of $\nnull(A(z))$ then it is easy to prove that $Q(z)=T(z)U(z)$ for some unimodular $U(z) \in \A^{p \times p}$. It follows that $\ker_\lambda A(z)$ is a subspace of $\C^n$. The next result generalizes \cite[Lemma 2.9]{DN}.

\begin{lemma}
$v \in \ker_\lambda A(z) \Leftrightarrow \exists w(z) \in \A^n$ : $A(z)w(z)=0$ and $w(\lambda)=v$.
\end{lemma}
\begin{proof}
Let $Q(z)$ be an invertible basis of $\nnull(A(z))$. Then from the existence of $w(z)$ we deduce that $w(z)=Q(z)c(z)$ for some $c(z) \in A^p$, and hence $v=Q(\lambda)c(\lambda)$. Conversely, if $v=Q(\lambda)c$ for some constant $c \in \C^p$ then set $w(z):=Q(z)c$.
\end{proof}

We are now in the position to extend some relevant definitions from \cite{DN,Nof12} to the case of analytic matrices. See also \cite{NV}.

\begin{definition}[Root vectors]
Given $A(z) \in \A^{m \times n}$ and $\lambda \in \Omega$, we say that $r(z) \in \A^n$ is a \emph{root vector} of order $k \geq 1$ at $\lambda$ for $A(z)$ if $r(\lambda) \not \in \ker_\lambda A(z)$ and $A(z)r(z)$ has a root of order $k$ at $\lambda$, i.e., $A(z)r(z)=(z-\lambda)^k v(z)$ for some $v(z) \in \A^m$, $v(\lambda) \neq 0$.
\end{definition}

\begin{definition}[$\lambda$-independent, complete, and maximal sets of root vectors]
Given $A(z) \in \A^{m \times n}$ and $\lambda \in \Omega$, suppose that $Q(z) \in \A^{n \times p}$ is an invertible basis of $\nnull(A(z))$ and that $\{r_i(z)\}_{i=1}^s \subset \A^n$ is a set of root vectors at $\lambda$ for $A(z)$ having orders $\{k_i\}_{i=1}^s$. Then:
\begin{itemize}
\item Such a set is $\lambda$-independent if $\begin{bmatrix}Q(\lambda) & r_1(\lambda) & \dots r_s(\lambda) \end{bmatrix} \in \C^{n \times (p+s)}$ has full column rank;
\item A $\lambda$-independent set is complete if there is no $\lambda$-independent set of strictly larger cardinality;
\item A complete set is ordered if $k_1 \geq \dots \geq k_s > 0$;
\item A complete ordered set is maximal if, for all $j$, there is no root vector $v(z)$ at $\lambda$ of order $k>k_j$ such that {$\begin{bmatrix}
Q(\lambda) & r_1(\lambda) & \dots & r_{j-1}(\lambda) & v(\lambda)
\end{bmatrix}$} has full column rank.
\end{itemize}
\end{definition}

At this point, we are well positioned to generalize several useful results found in \cite[Sections 3 and 4]{DN} and in \cite{DN2,NV} from polynomials (or rational functions) to analytic functions. We state these generalized results below as a number of propositions and theorems; in many cases, however, we omit the proofs as they are essentially the same as the proofs in \cite{DN,DN2,NV} for the polynomial (or rational) case. Here, ``essentially" points to minor modifications that, once made, allow us to follow the same main ideas of the original proofs. Typical examples of these minor modifications are: one may need to replace minimal bases with invertible bases, the ring of polynomials with the ring $\A$, the field of rational functions with the field of functions meromorphic on $\Omega$, or the expansion of a polynomial round a point to the Taylor series of an analytic function about a point in its domain of analyticity, etc. When, instead, a proof must be significantly different than its analogue appeared in \cite{DN,DN2,NV} for polynomial or rational matrices, we have included it.

\begin{theorem}[Generalization of Theorem 3.1 in \cite{DN}]\label{thm:31}
Let $S(z) \in \A^{m \times n}$ be in Smith form, and suppose that the rank of {$S(z)$} is $r$ and that the geometric multiplicity of $\lambda$ as an eigenvalue of $S(z)$ is $s$. Then $\{e_r,e_{r-1},\dots,e_{r-s+1}\}$, where $e_i \in \C^n$ is the $i$-th vector of the canonical basis, is a maximal set of root vectors at $\lambda$ for {$S(z)$}; moreover, their orders are the nonzero partial multiplicities of $\lambda$ as an eigenvalue of {$S(z)$}.
\end{theorem}

We note that \cite[Theorem 3.1]{DN} was stated referring to the local Smith form at $\lambda$ rather than to the (``global") Smith form as in Theorem \ref{thm:31}. The same could be done for analytic matrices, but for the sake of conciseness -- and since the results that follow can be equivalently deduced from either the local or global version of Theorem \ref{thm:31} -- we have decided here to state a global version rather than having to introduce the definition of local Smith form; interested readers can find it for instance in \cite{DN}.

\begin{lemma}[Generalization of Lemma 3.3 in \cite{DN}]
Let $P(z),Q(z) \in \A^{m \times n}$ and $\lambda \in \Omega$. Suppose that $Q(z)=A(z)P(z)B(z)$ for some $A(z) \in \A^{m \times m}, B(z) \in \A^{n \times n}$ and that $\det A(\lambda)\det B(\lambda)\neq 0$; moreover, let $M(z),N(z)$ be invertible bases for $\nnull(P(z)), \nnull(Q(z))$ respectively. Then $\ker_\lambda P(z)=\spann M(\lambda)=\spann B(\lambda)N(\lambda)$ and $\ker_\lambda Q(z)=\spann N(\lambda)=\spann \adj B(\lambda)M(\lambda)$.
\end{lemma}

{ 
\begin{proof}
It is clear that $\dim \ker_\lambda P(z) = \dim \ker_\lambda Q(z)$, because $A(z),B(z)$ are invertible over the field of meromorphic functions as otherwise $\det A(z) \det B(z)=0$, contradicting the assumption $\det A(\lambda) \det B(\lambda) \neq 0$. Moreover, we have $Q(z) \adj B(z) M(z) = \det B(z) P(z) M(z) = 0$ and thus the columns of $\adj B(z) M(z)$ lie in $\nnull(Q(z))$, implying that $\spann \adj B(\lambda) M(\lambda) \subseteq \ker_\lambda Q(z) = \spann N(\lambda)$. On the other hand, the constant matrix $\adj B(\lambda) = \det B(\lambda) B(\lambda)^{-1}$ is invertible, and hence $\rank \adj B(\lambda) M(\lambda) = \rank M(\lambda) = \rank N(\lambda)$. We conclude that $\adj B(\lambda) M(\lambda)$ is a basis for $\ker_\lambda Q(z)$.

An analogous argument can be used to show that $\ker_\lambda P(z)=\spann M(\lambda)=\spann B(\lambda)N(\lambda)$; we omit the details.
\end{proof}
  }

\begin{theorem}[Generalization of Theorem 3.4 in \cite{DN}]\label{thm:34}
Let $P(z),Q(z) \in \A^{m \times n}$ and $\lambda \in \Omega$. Suppose that $Q(z)=A(z)P(z)B(z)$ for some $A(z) \in \A^{m \times m}, B(z) \in \A^{n \times n}$ and that $\det A(\lambda)\det B(\lambda)\neq 0$. Then:
\begin{itemize}
\item If $\{v_i(z)\}_{i=1}^s$ are a maximal (resp. complete, $\lambda$-independent) set of root vectors at $\lambda$ for $Q(z)$ with orders $\ell_1 \geq \dots \ell_s > 0$, then {$\{B(z)v_i(z)\}_{i=1}^s$} are a maximal (resp. complete, $\lambda$-independent) set of root vectors at $\lambda$ for $P(z)$, with the same orders;
\item If $\{w_i(z)\}_{i=1}^s$ are a maximal (resp. complete, $\lambda$-independent) set of root vectors at $\lambda$ for $P(z)$ with orders $\ell_1 \geq \dots \ell_s > 0$, then {$\{\adj B(z)w_i(z)\}_{i=1}^s$} are a maximal (resp. complete, $\lambda$-independent) set of root vectors at $\lambda$ for $Q(z)$, with the same orders.
\end{itemize}
\end{theorem}

Note that Theorem \ref{thm:31} and Theorem \ref{thm:34} imply that every analytic matrix $A(z)$ has a maximal set of root vectors at $\lambda$, whose cardinality is precisely equal to the geometric multiplicity of $\lambda$ as an eigenvalue of $A(z)$. (To make full sense of the latter statement, we must formally agree that, if $\lambda$ is not an eigenvalue, then the empty set is a set of root vectors at $\lambda$ for $A(z)$, and in fact the only such possible set.)

\begin{theorem}[Generalization of Theorem 3.10 in \cite{NV}]\label{thm:310}
Let $A(z) \in \A^{m \times n}$ and suppose that $Q(z) \in \A^{n \times p}$ is an invertible basis for $\nnull(A(z))$. Suppose that $\{v_i(z)\}_{i=1}^s$ are root vectors at $\lambda \in \Omega$ for $A(z)$. Then, they are a complete set if and only if the set of the columns of the matrix 
\[ B : = \begin{bmatrix}
Q(\lambda) & v_1(\lambda) & \dots v_s(\lambda) 
\end{bmatrix} \]
is  a basis for $\ker A(\lambda)$.
\end{theorem}
\begin{proof}
By definition, $A(z)Q(z)=0$ and $A(z)v_i(z)=(z-\lambda)^{k_i} w_i(z)$, where $k_i$ is the order of $v_i(z)$ as a root vector of $A(z)$, implying $A(\lambda)B=0$. If $\{v_i(z)\}_{i=1}^s$ is complete, then it is in particular $\lambda$-independent so $B$ has full column rank. If the set of the columns of $B$ is not a a basis for $\ker A(\lambda)$, we can complete it to a basis by adding some (constant) vectors $\{ x_i \}_{i=1}^c$. But then $x_i$ are root vectors at $\lambda$ for $A(z)$, and by construction $\{v_i(z)\}_{i=1}^s \cup \{x_i\}_{i=1}^c$ are a $\lambda$-independent set, contradicting completeness of the original set. Conversely, if the set of root vectors is not complete then the set of the columns of $B$ cannot be a basis of $\ker A(\lambda)$ for dimensional reasons, as by starting from a complete set of root vectors we can construct $t > p + s$ linearly independent vectors in $\ker A(\lambda)$.
\end{proof}

\begin{theorem}[Generalization of Theorem 4.1 and Theorem 4.2 in \cite{DN}]\label{thm:410}
Suppose that $A(z) \in \A^{m \times n}$ has partial multiplicities $m_1 \geq \dots \geq m_s$ at $\lambda \in \Omega$. Then:
\begin{enumerate}
\item All complete sets of root vectors of $A(z)$ at $\lambda$ have the same cardinality, equal to $s$;
\item All maximal sets of root vectors of $A(z)$ at $\lambda$ have the same ordered list of orders $m_1 \geq \dots \geq m_s$, equal in turn to the partial multiplicities of $\lambda$ as an eigenvalue of $A(z)$;
\item If $\{v_i(z)\}_{i=1}^s$ is an ordered complete set of root vectors at $\lambda$ for $A(z)$ having orders $\ell_1 \geq \dots \geq \ell_s$, then
\begin{itemize}
\item[3.1] $m_i \geq \ell_i$ for all $i=1,\dots,s$;
\item[3.2] $\{v_i(z)\}_{i=1}^s$ is a maximal set of root vectors at $\lambda$ for $A(z)$ if and only if $\ell_i=m_i$ for all $i=1,\dots,s$;
\item[3.3] $\{v_i(z)\}_{i=1}^s$ is a maximal set of root vectors at $\lambda$ for $A(z)$ if and only if $\sum_{i=1}^s m_i = \sum_{i=1}^s \ell_i$.
\end{itemize}
\end{enumerate}
\end{theorem}

\begin{theorem}[Generalization of Lemma 5.2 in \cite{DN2}]
Suppose that $A(z) \in \A^{m \times n}$ has partial multiplicities $m_1 \geq \dots \geq m_s$ at $\lambda \in \Omega$. Moreover, let $\{v_i(z)\}_{i=1}^c$ be a $\lambda$-independent set of root vectors at $\lambda$ for $A(z)$ having orders $\ell_1 \geq \dots \geq \ell_c$. Then $c \leq s$ and $\ell_i \leq m_i$ for all $i=1,\dots,c$.
\end{theorem}

The concept of root vectors allows us also to give appropriate definitions of eigenvectors, even for rectangular or square but singular analytic matrices. To this goal, we give below some relevant definitions that generalize analogous ones given (for polynomial or rational matrices) in \cite{DN,NV}. Fix an analytic matrix $A(z) \in \A^{m \times n}$ and a scalar $\lambda \in \Omega$, and note that generally $\ker_\lambda A(z) \subseteq \ker A(\lambda)$. However, the inclusion is strict if and only if $\lambda$ is an eigenvalue of $A(z)$ if and only if the cardinality of any maximal set of root vectors at $\lambda$ for $A(z)$ is positive. This motivates the following definition.

\begin{definition}\label{def:evec}
An equivalence class $[v] \in \ker A(\lambda)/\ker_\lambda A(z)$ is called a right eigenvector of $A(z)$ associated with $\lambda$ if $[v]\neq[0]$.
\end{definition}

Following the same arguments as in \cite{DN,NV}, we see that Definition \ref{def:evec} implies that eigenvectors associated with $\lambda$ exists if and only if $\lambda$ is an eigenvalue. Theorem \ref{thm:310} implies moreover that $[v]$ is an eigenvector if and only if
\[ [0] \neq [v] \in \spann \{ [x_1(\lambda)],\dots,[x_s(\lambda)] \} \]
where $\{x_i(z)\}_{i=1}^s$ is any complete set of root vectors at $\lambda$ for $A(z)$. Equivalently, $[v]$ is an eigenvector if and only if $v=\sum_{i=1}^s x_i(\lambda) c_i + w$ for some, { not all zero,} coefficients $c_i \in \C$ and some vector $w \in \ker_\lambda A(z)$. When $\lambda$ is a simple eigenvalue, it is possible to make $v$ unique up to a phase, by asking that $v \in \ker_\lambda A(z)^\perp$ and that $\|v\|_2=1$; in the polynomial case, this was useful for example in \cite{KG,LN} in the context of conditioning analysis and in \cite{HMP,KG} to devise practical algorithms. We expect that similar developments may be possible also in the analytic case. We conclude by noting that, if $A(z)$ has full column rank, then $\ker_\lambda A(z) = \{0\}$ and thus we recover the familiar definition of an eigenvector being a nonzero vector $v$ such that $A(\lambda)v=0$.

\begin{remark}
In this paper, we focused on extending the theory of root vectors from polynomial \cite{DN} to analytic matrices, by using invertible bases and the theory of modules. In \cite{NV}, the theory of root vectors was extended from polynomial to rational matrices by using the theory of discrete valuations. Adding valuation theory to the tools of this paper makes it similarly possible to state a theory of root vectors for meromorphic matrices. We opted to not do it in this paper for two reasons: (1) because meromoprhic matrices do not appear to be as central, in the current linear algebra literature, as polynomial, analytic or rational matrices (2) for simplicity of exposition, i.e., to avoid an excessive density of tools that are possibly not too familiar to every matrix theorist or numerical linear algebraist. We thus leave the task to future research. 
\end{remark}

\subsection{Example}

In this subsection, we work out an example to illustrate the results of Section \ref{sec:maximalsets}. Let $\Omega=\C$ so that $\A$ is the ring of entire functions, and consider the analytic matrix
\[  A(z) = \begin{bmatrix}
2 z e^{2z} & z e^z & z \sinh(z)\\
e^{2z} (\sin(z)^2-2z) & e^z(\sin(z)^2-z) & e^z \cos(z) \sin(z)^2-z \sinh(z)\\
e^z(2ze^z-\sin(z)^2) & ze^z-\sin(z)^2 & \cos(z)^3-\cos(z)+z\sinh(z)
\end{bmatrix}.\]
We are going to study the spectral properties of $A(z)$ at its eigenvalue $\lambda=0$: in particular we will exhibit a maximal set of root vectors, thus obtaining the partial multiplicities of $0$ as well as corresponding eigenvectors  (according to Definition \ref{def:evec}). As a preliminary step, note that a straightforward computation shows that $\det A(z)=0$, and that the leading principal $2 \times 2$ minor of $A(z)$ is $e^{3z} z \sin(z)^2 \neq 0$. Hence,  $\rank A(z)=2$. It can be readily verified that an invertible basis of $\nnull(A(z))$ is
\[ Q(z)=\begin{bmatrix}
e^z\cos(z)-\sinh(z)\\
e^z(\sinh(z)-2e^{z}\cos(z))\\
e^{2z}
\end{bmatrix}  \Rightarrow \ker_0 A(z)=\mathrm{span} \left( \begin{bmatrix}
1\\-2\\1
\end{bmatrix} \right) .\]
Let now
\[ v_1(z)=\begin{bmatrix}
1\\
-e^z\\
0
\end{bmatrix}, \qquad v_2(z)=\begin{bmatrix}
1\\
-2 e^z\\
0
\end{bmatrix}.\]
We claim that $\{v_1(z),v_2(z)\}$ is a maximal set of root vectors at $0$ for $A(z)$. Note first that
\[ A(z)v_1(z)=ze^{2z} \begin{bmatrix}
1\\-1\\1
\end{bmatrix}, \qquad A(z)v_2(z)=e^z \sin(z)^2\begin{bmatrix}
0\\
-e^z\\
1
\end{bmatrix},\]
and that
\[ v_1(0)=\begin{bmatrix}
1\\-1\\0
\end{bmatrix} \not\in \ker_0 A(z), \qquad v_2(0)=\begin{bmatrix}
1\\-2\\0
\end{bmatrix} \not\in \ker_0 A(z),\]
and hence $v_1(z),v_2(z)$ are root vectors at $0$ for $A(z)$ of order $1$ and $2$ respectively: note indeed the Taylor expansions $ze^z=z+o(z)$ and $e^z \sin(z)^2 = z^2 + o(z)$. The set $\{v_i(z)\}_{i=1}^2$ is $0$-independent because
\[ \mathrm{rank} \begin{bmatrix}
Q(0) & v_1(0) & v_2(0)
\end{bmatrix} = \mathrm{rank} \begin{bmatrix}
1&1&1\\
-2&-1&-2\\
1&0&0
\end{bmatrix}=3.\]
Moroever, we have $A(0)=0$ and hence, by Theorem \ref{thm:310}, the set $\{v_i(z)\}_{i=1}^2$ is complete. Finally, one can compute a Smith form of $A(z)$, for example by computing determinantal divisors (GCDs of minors of a given size); indeed, recall from the introduction that the invariant factors are the ratios of the determinantal divisor. We thus obtain that the invariant factors of $A(z)$ are $z,\sin(z)^2,0$. This implies that the partial multiplicities of $0$ are $1,2$, and by Theorem \ref{thm:410} the set $\{v_i(z)\}_{i=1}^2$ is therefore maximal. Moreover, $[v_1(0)]$ and $[v_2(0)]$, seen as equivalence classes $\mod \ker_0 A(z)$, are linearly independent eigenvectors of $A(z)$ associated with the eigenvalue $0$.
\section*{Acknowledgements}

I thank Froil\'{a}n Dopico for reading a preliminary version of this paper and sharing very useful remarks. {I am also very grateful to an anonymous reviewer for many careful suggestions. In particular, I am indebted to the referee for Remarks \ref{rem:referee1} and \ref{rem:referee2} as well as for other improvements of the presentation.}

\end{document}